\documentclass{IEEEtran}
\usepackage{cite}
\usepackage{amsmath,amssymb,amsfonts}
\usepackage{algorithmic}
\usepackage{graphicx}
\usepackage{textcomp}
\def\BibTeX{{\rm B\kern-.05em{\sc i\kern-.025em b}\kern-.08em
    T\kern-.1667em\lower.7ex\hbox{E}\kern-.125emX}}

\usepackage{enumerate}
\usepackage{flushend}

\newcommand{\set}[1]{\left\{#1\right\}}
\newcommand{\norm}[1]{\left\Vert #1 \right\Vert}

\newcommand{\ra}{\rightarrow}
\newcommand{\Real}{\mathbb{R}}
\newcommand{\eps}{\varepsilon}
\newcommand{\B}{\mathbb{B}}
\renewcommand{\S}{\mathcal{S}}
\renewcommand{\subset}{\subseteq}
\newcommand{\Sd}{\mathcal{S}_{\delta}}
\newcommand{\dtS}{\mathcal{DTS}}
\newcommand{\dtSd}{\mathcal{DTS}_{\delta}}
\newcommand{\K}{\mathcal{K}}

\newtheorem{thm}{Theorem}

\newtheorem{lem}[thm]{Lemma}

\newtheorem{assum}[thm]{Assumption}
\newtheorem{prop}[thm]{Proposition}

\newtheorem{defn}[thm]{Definition}
\newtheorem{exmp}[thm]{Example}
\newtheorem{rem}[thm]{Remark}

\begin{document}
\title{Converse Barrier Functions via Lyapunov Functions}

\author{Jun Liu, \IEEEmembership{Senior Member, IEEE}
\thanks{This work was supported in part by the Natural Sciences and Engineering Research Council of Canada, the Canada Research Chairs program, and the Government of Ontario through an Early Researcher Award.}
\thanks{Jun Liu is with the Department of Applied Mathematics, University of Waterloo, Waterloo, Ontario N2L 3G1, Canada (e-mail: j.liu@uwaterloo.ca). }
}

\maketitle

\begin{abstract}
We prove a robust converse barrier function theorem via the converse Lyapunov theory. While the use of a Lyapunov function as a barrier function is straightforward, the existence of a converse Lyapunov function as a barrier function for a given safety set is not. We establish this link by a robustness argument. We show that the closure of the forward reachable set of a robustly safe set must be robustly asymptotically stable under mild technical assumptions. As a result, all robustly safe dynamical systems must admit a robust barrier function in the form of a Lyapunov function for set stability. We present the results in both continuous-time and discrete-time settings and remark on connections with various barrier function conditions.

\end{abstract}


\section{Introduction}
\label{sec:introduction}

The use of barrier functions to ensure set invariance and safety in control of dynamical systems has gained popularity in recent years in safety-critical control applications \cite{prajna2007convex,prajna2005necessity,prajna2004safety,prajna2007framework,wieland2007constructive,agrawal2017discrete,cheng2019end,ahmadi2019safe,ames2016control,xu2015robustness,tee2009barrier,romdlony2016stabilization}. The readers are referred to \cite{ames2016control} for a nice introduction on the background of barrier functions. 

From the earlier work \cite{prajna2005necessity,prajna2007convex} to recent results \cite{wisniewski2016converse,ratschan2018converse}, converse theorems for barrier functions played an important role in understanding how safety properties can indeed be characterized by barrier functions. The more stringent conditions in  \cite{prajna2005necessity,prajna2007convex} for the existence of converse barrier functions are relaxed in \cite{wisniewski2016converse} to a class of structurally table dynamical systems (more precisely, Morse-Smale vector fields) and in \cite{ratschan2018converse} to a robust safety requirement. 

In this paper, inspired by the recent work \cite{ratschan2018converse} and the connections made in \cite{xu2015robustness} (see also \cite{ames2016control}) between a barrier function and a Lyapunov function, we prove that, for all \textit{robustly} safe dynamical systems, barrier functions can be constructed from Lyapunov functions. The use of Lyapunov functions to ensure set invariance is standard  \cite{khalil2002nonlinear} (see also \cite{prajna2007convex}). The authors of \cite{xu2015robustness,ames2016control} also highlighted that if the barrier function conditions are satisfied in a neighborhood of the safety set, then the barrier function can indeed be regarded as a Lyapunov function. What is missing, however, are conditions under which such barrier functions exist assuming safety of the system. We establish this link by proving that the closure of the robust reachable set of a robustly safe set must be robustly asymptotically stable under mild technical assumptions (Theorem \ref{thm:stable}). The results of this paper could help provide a potentially more unified view of the Lyapunov function and barrier function theories, because how to simultaneously satisfy Lyapunov and barrier function conditions are important in practice and but technically challenging \cite{ames2016control,romdlony2016stabilization}. 

\textbf{Notation:} For $x\in\Real^n$ and $r\ge 0$, we denote the ball of radius $r$ centered at $x$ by $B_{r}(x)=\set{y\in\Real^n:\,\norm{y-x}\le r}$, where $\norm{\cdot}$ is the Euclidean norm. For a closed set $A\subset\Real^n$ and $x\in\Real^n$, we denote the distance from $x$ to $A$ by $\norm{x}_{A}=\inf_{y\in A}\norm{x-y}$ and $r$-neighborhood of $A$ by $B_r(A)=\cup_{x\in A}B_r(x)=\set{x\in\Real^n:\,\norm{x}_A\le r}$. For convenience, we also write $\B=B_1(0)$ and $r\B=B_r(0)$.

The remainder the paper is organized as follows. We present some preliminaries on barrier and Lyapunov functions for continuous-time systems in Section \ref{sec:prel}. We prove a converse barrier function theorem by a converse Lyapunov function theorem in Section  \ref{sec:main}. The results of Section \ref{sec:main} are extended to discrete-time systems in Section \ref{sec:discrete}. The paper is concluded in Section \ref{sec:conclusion}.

\section{Preliminaries}\label{sec:prel}

Consider a continuous-time dynamical system 
\begin{equation}\label{eq:sys}
    x' = f(x),
\end{equation}
where $x\in\Real^n$ and $f:\,\Real^n\ra\Real^n$ is assumed to be locally Lipschitz. %
For each $x_0\in\Real^n$, we denote the unique solution starting from $x(0)=x_0$ and  defined on the maximal interval of existence by $x(t;x_0)$ or simply $x(t)$ if $x_0$ is not emphasized.

Given a scalar $\delta\ge 0$, a \emph{$\delta$-perturbation} of the dynamical system (\ref{eq:sys}) is described by the differential inclusion
\begin{equation}\label{eq:p1}
    x' \in F_{\delta}(x),
\end{equation}
where $F_{\delta}(x)=B_{\delta}(f(x))$. An equivalent description of the  $\delta$-perturbation of system (\ref{eq:sys}) can be given by
\begin{equation}\label{eq:p2}
    x'(t) = f(x(t)) + d(t),
\end{equation}
where $d:\,\Real\ra\delta\B$ is any measurable signal. We denote system (\ref{eq:sys}) by $\S$ and its $\delta$-perturbation by $\Sd$. Note that $\Sd$ reduces to $\S$ when $\delta=0$. A solution of $\Sd$ starting from $x(0)=x_0$ can be denoted by $x(t;x_0,d)$ or simply $x(t)$, where $d$ is a given disturbance signal. The set of all solutions for $\Sd$ starting from $x_0$ is denoted by $S_{\delta}(x_0)$. We are only interested in forward solutions (i.e., solutions defined in positive time) in this paper. Set invariance, defined below and used in this paper, also only concerns forward invariance. 

\begin{defn}[Invariant set]
A set $\Omega\subset\Real^n$ is said to be an \emph{invariant set} of $\Sd$ if all solutions of $\Sd$ starting in $\Omega$ remain in $\Omega$ in positive time.  
\end{defn}

\begin{defn}[Robustly invariant set]
A set $\Omega\subset\Real^n$ is said to be a \emph{$\delta$-robustly invariant set} of $\S$ if it is an invariant set of $\Sd$ for some $\delta\ge 0$. It is said to be a \emph{robustly invariant set} of $\S$ if it is a $\delta$-robustly invariant set for some $\delta>0$. 
\end{defn}

\begin{defn}[Robustly safe set]
Given an \textit{unsafe} set $U\subseteq\Real^n$, a set $W\subseteq\Real^n$ is said to be \textit{$\delta$-robustly safe} w.r.t. to $U$ if all solutions of $\Sd$ starting from $W$ will not enter $U$. 
\end{defn}

An immediate connection between robustly safe and invariant sets is the following. 
\begin{prop}\label{prop:invariance}
If there exists a $\delta$-robustly invariant set $\Omega$ such that $W\subset\Omega$ and $\Omega\cap U=\emptyset$, then $W$ is $\delta$-robustly safe w.r.t. to $U$. 
\end{prop}

\begin{defn}[Robust barrier function]\label{def:barrier}
Given sets $W,U\subset\Real^n$, a continuously differentiable function  $B:\,\Real^n\ra\Real$ is said to be a \textit{$\delta$-robust barrier function} for $W$ and $U$ if the following conditions are satisfied: 
\begin{enumerate}[(1)]
    \item $B(x)\ge 0$ for all $x\in W$;
    \item $B(x)<0$ for all $U$; and 
    \item $\nabla B(x)\cdot(f(x)+d)>0$ for all $x$ such that $B(x)=0$ and all $d\in\delta\B$. 
\end{enumerate}
\end{defn}

\begin{rem}
The choice of sign for $B$ to indicate a safe set is rather arbitrary because we can always negate it. Here we use the condition $B(x)\ge 0$ to describe the safe set (the same as \cite{ames2016control} and \cite{ratschan2018converse}) instead of $B(x)\le 0$ (as in the original work \cite{prajna2004safety}). 
\end{rem}

A $\delta$-robust barrier function for $(W,U)$ provides a certificate for $\delta$-robust safety of $W$ w.r.t. $U$, as summarized in the following result. 

\begin{prop}[Sufficiency of barrier functions \cite{prajna2004safety,prajna2005necessity}]
If there exists a $\delta$-robust barrier function for $(W,U)$, then $W$ is $\delta$-robustly safe w.r.t. $U$. 
\end{prop}

A proof that leads to a slightly different conclusion can be found in \cite{prajna2005necessity}. We provide a short proof below for completeness. 

\begin{IEEEproof}
We show that the set $C=\set{x\in\Real^n:\,B(x)\ge 0}$ is an invariant set for $\Sd$. Robust safety of $W$ follows immediately in view of Proposition \ref{prop:invariance} and conditions (1)--(2) of Definition \ref{def:barrier}. Suppose that $C$ is not invariant for $\Sd$. Then there exists a solution $x(\cdot)$ for $\Sd$ such that $x(0)\in C$ and $x(t)\not\in C$ for some $t>0$. Define 
$$
\overline{t}=\sup\{t\ge 0:\,x(t)\in C\}.
$$
Then $\overline{t}$ is well defined and finite. By continuity of $B(x(t))$, we have $B(x(\overline{t}))=0$. This implies that 
$$\frac{dB(x(t))}{dt}=\nabla B(x(t))\cdot (f(x(t)) + d(t))>0$$ 
at $t=\overline{t}$. Hence, for $\eps>0$ sufficiently small, we have $B(x(t))>B(x(\overline{t}))=0$ for $t\in (\overline{t},\overline{t}+\eps]$. This contradicts the definition of $\overline{t}$. 
\end{IEEEproof}

\begin{rem}
Note that the strict inequality $\nabla B(x)\cdot(f(x)+d)>0$ is needed to guarantee the set $\set{x\in\Real^n:\,B(x)\ge 0}$ is forward invariant. The original paper  \cite{prajna2004safety} used the non-strict inequality condition: $\nabla B(x)\cdot(f(x)+d)\ge 0$ for all $x$ such that $B(x)=0$. This condition has been known to be unsound (see, e.g., \cite[Example 2]{taly2009deductive}; see also \cite[Remark after Theorem 3]{gard1980strongly}). Safety properties of dynamical systems are intimately related to set invariance, on which there is a rich history of investigation (interested readers can refer to \cite{blanchini1999set} for more information; see also \cite[Section 3]{gard1980strongly}). 
\end{rem}

Several converse theorems for barrier functions have been proved in the literature \cite{prajna2005necessity,prajna2007convex,wisniewski2016converse,ratschan2018converse}. We quote a most recent result by Ratschan as follows. 

\begin{thm}[Necessity of barrier functions \cite{ratschan2018converse}]\label{thm:rats}
Suppose that the closure of $W$ and $U$ are disjoint and the complement of $U$ is bounded. If $W$ is $\delta$-robustly safe w.r.t. $U$, then there exists a continuously differentiable function  $B:\,\Real^n\ra\Real$ satisfying the following conditions: 
\begin{enumerate}[(1)]
    \item $B(x)\ge 0$ for all $x\in W$;
    \item $B(x)<0$ for all $U$; and 
    \item $\nabla B(x)\cdot f(x)>0$ for all $x$ such that $B(x)=0$. 
\end{enumerate}
\end{thm}

While condition (3) appears to be slightly different from item (3) in Definition \ref{def:barrier}, we will remark on the connections between them, as well as with other variants of barrier function conditions, in Section \ref{sec:main} (see Remark \ref{rem:barrier}). 

We say a continuous function $\alpha:[0,a)\ra\Real$ belongs to class $\K$ and write $\alpha\in\K$ if $\alpha$ is strictly increasing and $\alpha(0)=0$. %

\begin{defn}[Set stability]\label{def:stability}
A closed set $A\subset\Real^n$ is said to be \textit{$\delta$-robustly uniformly asymptotically stable} ($\delta$-RUAS) for $\S$ if the following two conditions are met:
\begin{enumerate}[(1)]
    \item For every $\eps>0$, there exists a $\delta_{\eps}>0$ such that $\norm{x(0)}_A<\delta_\eps$ implies $\norm{x(t)}_A<\eps$ for any solution $x(t)$ of $\Sd$; and 
    \item There exists some $\rho>0$ such that, for every $\eps>0$, there exists some $T>0$ such that $\norm{x(t)}_A<\eps$ for any solution $x(t)$ of $\Sd$ whenever $\norm{x(0)}_A<\rho$ and $t\ge T$. 
\end{enumerate}
\end{defn}

It is not difficult to see that a $\delta$-robustly uniformly asymptotically stable set $A$ must be $\delta$-robustly invariant. 

\begin{defn}[Robust Lyapunov function]\label{def:lyap}
Let $D\subset\Real^n$ be an open set containing a closed set $A\subset\Real^n$. A continuously differentiable function $V:\,D\ra\Real$ is said to be a \textit{$\delta$-robust Lyapunov function} for $\S$ w.r.t. $A$ if the following two conditions are satisfied:
\begin{enumerate}[(1)]
    \item there exist class $\K$ functions $\alpha_1$ and $\alpha_2$ such that 
    $$
    \alpha_1(\norm{x}_A) \le V(x) \le \alpha_2(\norm{x}_A)
    $$
    for all $x\in D$; and 
    \item there exists a class $\K$ functions $\alpha_3$ such that
    $$
    \nabla V(x)\cdot((f(x)+d)\le -\alpha_3(\norm{x}_A)
    $$
    for all $x\in D$ and $d\in \delta\B$. 
\end{enumerate}
\end{defn}

There are well-known Lyapunov characterizations of set stability. 

\begin{thm}[Lyapunov characterization of set\label{thm:lyap} stability\cite{wilson1969smoothing,lin1996smooth}]
A closed set $A\subset\Real^n$ is $\delta$-RUAS for $\S$ if and only if there exists a $\delta$-robust Lyapunov function for $\S$ w.r.t. $A$. 
\end{thm}

\section{Robust converse barrier functions via Lyapunov functions}
\label{sec:main}

In this section, we prove a version of converse barrier function theorem by resorting to converse Lyapunov theory. 

We first introduce some notation. Let $R_{\delta}^t(x_0)$ denote the set reached by solutions of $\Sd$ at time $t$ starting from $x_0$, i.e.,
$$
R_{\delta}^t(x_0) = \set{x(t):\,x(\cdot)\in S_{\delta}(x_0))}. 
$$
We further define 
$$
R_{\delta}(x_0) = \bigcup_{t\ge 0} R_{\delta}^t(x_0),
$$
and, for a set $W\subset\Real^n$, 
$$
R_{\delta}^t(W) =  \bigcup_{x_0\in W} R_{\delta}^t(x_0),\quad R_{\delta}(W) =  \bigcup_{x_0\in W} R_{\delta}(x_0). 
$$
Clearly, $R_\delta(W)$ is a $\delta$-robustly invariant set of $\S$. If $W$ is $\delta$-robustly safe, then $R_\delta(W)\cap U\neq\emptyset$. If the complement of $U$ is bounded (as assumed in Theorem \ref{thm:rats}), then $R_\delta(W)$ is bounded. Let $\Omega=\overline{R_\delta(W)}$. Then $\Omega$ is compact. Without further assumption, $\Omega$ may intersect with $U$ as shown in the following example. 

\begin{exmp}\label{ex:counter}
Consider $\S$ defined by $x' = - x$. Let $W=[-0.1,0.1]$ and $\delta=0.2$. Then $R_\delta(W)=(-0.2,0.2)$ and $\Omega=[-0.2,0.2]$. If $U=(-\infty,-2]\cup[2,\infty)$, then $W$ is $\delta$-robustly safe w.r.t. $U$ because $R_\delta(W)\cap U=\emptyset$. Yet $\Omega\cap U\neq\emptyset$. 
\end{exmp}

Note that the assumptions of Theorem \ref{thm:rats} are indeed satisfied by the example above. While additionally assuming $U$ to be open  will lead to $\Omega\cap U=\emptyset$, we need a slightly stronger assumption for the purpose of this section, that is, $\Omega\cap \overline{U}=\emptyset$. This is summarized in the following assumption. 

\begin{assum}\label{as:main}
The set $W$ is $\delta$-robustly safe w.r.t. $U$ and $\Omega\cap \overline{U}=\emptyset$, where $\Omega=\overline{R_\delta(W)}$. Furthermore, either $\Omega$ is compact or $f$ is globally Lipschitz. 
\end{assum}

With this assumption, we prove the following result on converse barrier functions.

\begin{thm}[Robustly safe sets admit robust barrier functions]\label{thm:main}
Suppose that Assumption \ref{as:main} holds. Then for any $\delta'\in(0,\delta)$, there exists a $\delta'$-robust barrier function for $(W,U)$.
\end{thm}

The conclusion of the above result is slightly stronger than the main result in \cite{ratschan2018converse} (quoted as Theorem \ref{thm:rats} in Section \ref{sec:prel} above) in two aspects: (1) we show the existence of a $\delta'$-robust barrier function for any $\delta'\in(0,\delta)$; (2) we do not assume $\Omega$ to be compact, when $f$ is globally Lipschitz\footnote{In fact, $f$ being Lipschitz in a neighborhood of $\Omega$ suffices.}. Assumption \ref{as:main} appears to be stronger than that of Theorem \ref{thm:rats} in that it requires $\Omega\cap \overline{U}=\emptyset$. Nonetheless, the proof of Theorem \ref{thm:rats} (see, e.g., Lemma 5 in \cite{ratschan2018converse}) seems to be using this fact without explicitly mentioning or proving it. Example \ref{ex:counter} above shows that this does not readily follow from the assumptions of Theorem \ref{thm:rats}. Despite these subtle technical differences, the main message of this section, however, is that converse barrier functions can be constructed from Lyapunov functions. 

The construction relies on showing that the closure of the reachable set of the robustly safe set, i.e., the set $\Omega=\overline{R_\delta(W)}$, is robustly asymptotically stable (Theorem \ref{thm:stable}). The following technical lemma on reachable sets of a perturbed system plays an important role in proving Theorem \ref{thm:stable}.

\begin{lem}\label{lem:control}
Fix any $\delta'\in(0,\delta)$ and $\tau>0$. Let $K\subset\Real^n$ be a compact set. Then there exists some $r=r(K,\tau,\delta',\delta)>0$ such that the following holds: if there is a solution $x$ of $\S_{\delta'}$ such that $x(s)\in K$ for all $s\in[0,T]$, where $T\ge\tau$, then for any $y_0\in B_r(x(0))$ and any $y_1\in B_r(x(T))$, we have $y_1\in R_{\delta}^T(y_0)$, i.e., $y_1$ is reachable at $T$ from $y_0$ by a solution of $\Sd$. Furthermore, if $f$ is globally Lipschitz, $r$ can be chosen to be independent of $K$.
\end{lem}

\begin{IEEEproof}
Let 
$$
y(s) = x(s) + \frac{s}{T}[y_1 - x_1 + (x_0-y_0)] + (y_0-x_0),\quad s\in[0,T].
$$
Then $y(0)=y_0$ and $y(T)=y_1$. Furthermore, 
\begin{align*}
    \norm{y(s)-x(s)} &\le \norm{y_1-x_1}\frac{s}{T} + \norm{y_0-x_0}(1-\frac{s}{T}) \\ 
    &\le r (\frac{s}{T} + 1  - \frac{s}{T}) = r,
\end{align*}
and
\begin{align*}
    \norm{y'(s)-x'(s)} &\le \norm{\frac{1}{T}[y_1 - x_1 + (x_0-y_0)] } \\
    & \le \frac{1}{T}[\norm{y_1-x_1}+\norm{x_0-y_0}] \le \frac{2r}{T},
\end{align*}
for all $s\in [0,T]$. Hence
\begin{align*}
    &\norm{y'(s)-f(y(s))} \\
    &= \norm{y'(s)-x'(s) + x'(s) - f(x(s)) + f(x(s)) - f(y(s))} \\
    & \le \frac{2r}{T} + \delta' + L r,
\end{align*}
where we used the triangle inequality, the fact that $x$ is a solution of $\S_{\delta'}$, and Lipschitz continuity of $f$ on the set $B_r(K)$. By picking $r$ sufficiently small such that $\frac{2r}{\tau} + \delta' + L r<\delta$, then we have
$\norm{y'(s)-f(y(s))}<\delta$ for all $s\in[0,T]$. Thus $y$ is a solution of $\Sd$ and the conclusion follows. Note that the choice of $r$ only depends on $K$, $\tau$, $\delta'$, and $\delta$. The dependence on $K$ is removed if $f$ is globally Lipschitz.
\end{IEEEproof} 

\begin{rem}
Lemma \ref{lem:control} extends the statement of Lemma 1 in \cite{ratschan2018converse}, where the proof was omitted. Lemma \ref{lem:control} is slightly stronger because it says that we can steer any point in a small neighborhood of $x(0)$ (as opposed to only $x(0)$) to a small neighborhood of $x(T)$. This fact is needed in the proof of Theorem \ref{thm:stable} below. Lemma \ref{lem:control} also allows $T$ to vary as long as it is lower bounded by $\tau$. The proof given here is elementary and constructive. Similar argument (of a simpler version) appeared in the proof of Theorem 1 in \cite{liu2019robust}. 
\end{rem}

\begin{thm}[Robustly invariant sets are robustly asymptotically stable]\label{thm:stable}
If Assumption \ref{as:main} holds, then for any $\delta'\in[0,\delta)$, the set $\Omega=\overline{R_\delta(W)}$ is $\delta'$-RUAS for $\S$. 
\end{thm}

\begin{IEEEproof}
We verify conditions (1) and (2) of Definition \ref{def:stability}. 

(1) For any $\eps>0$, let $\tau>0$ be the minimal time that is required for solutions of $\S_{\delta'}$ to escape from $B_{\frac{\eps}{2}}(\Omega)$ to $B_{\eps}(\Omega)$. The existence of such a $\tau$ follows from that $f$ is locally Lipschitz and an argument using Gronwall's inequality. Pick  $\delta_\eps<\min(r,\frac{\eps}{2})$, where $r$ is from Lemma \ref{lem:control}, applied to the set $B_{\eps}(\Omega)$ and scalars $\tau$, $\delta'$, and $\delta$. Let $x$ be any solution of $\S_{\delta'}$ such that $\norm{x(0)}_\Omega<\delta_\eps$. We show that $\norm{x(t)}_\Omega<\eps$ for all $t\ge 0$. Suppose that this is not the case. Then $\norm{x(t_1)}_\Omega\ge\eps$ for some $t_1\ge\tau>0$.  
Since $\delta_\eps<r$, we can always pick $y_0\in R_\delta(W)$ such that $y_0\in B_{r}(x(0))$ by the triangle inequality. %
By Lemma \ref{lem:control}, there exists a solution of $\Sd$ from $y_0\in R_\delta(W)$ to $x(t_1)\not\in\Omega$. This contradicts that $R_\delta(W)$ is $\delta$-robustly invariant. 

(2) Fix any $\eps_0>0$. Following part (1), choose $\delta_{\eps_0}$ such that $\norm{x(0)}_\Omega<\delta_{\eps_0}$ implies $\norm{x(t)}_\Omega<\eps_0$ for any solution $x(t)$ of $\S_{\delta'}$. Let $r$ be chosen according to Lemma \ref{lem:control} with the set $B_{\eps_0}(\Omega)$ and scalars $\tau=1$, $\delta'$, and $\delta$. Choose $\rho\in(0,r)$. Let $x$ be any solution of $\S_{\delta'}$. We show that $\norm{x(0)}_\Omega<\rho$ implies that $x(t)\in R_\delta(W)$ for all $t\ge 1$. Suppose that this is not the case. Then there exists some $t_1\ge 1$ such that $x(t_1)\in\partial \Omega$ or $x(t_1)\not\in\Omega$. In either case, we can pick $y_1\in B_r(x(t_1))$ such that $y_1\not\in\Omega$ and $y_0\in B_r(x(0))$ such that $y_0\in  R_\delta(W)$. By Lemma \ref{lem:control}, there exists a solution of $\Sd$ from $y_0\in R_\delta(W)$ to $y_1\not\in\Omega$. This contradicts that $R_\delta(W)$ is $\delta$-robustly invariant. Hence $x(t)\in R_\delta(W)\subset\Omega$ for all $t\ge 1$. %
This clearly implies (2). 
\end{IEEEproof}

The conclusion of Theorem \ref{thm:stable} cannot be strengthened in the sense that the set $\Omega=\overline{R_\delta(W)}$ may not be $\delta$-RUAS for $\S$ as shown in the simple example below.

\begin{exmp}
Consider $\S$ defined by $x' = - x + x^2$. Let $W=[-0.1,0.1]$ and $\delta=0.25$. Then $R_\delta(W)=(-0.25,0.25)$ and $\Omega=[-0.25,0.25]$. Solutions of $\Sd$ starting from $x_0=0.25+\eps$, where $\eps>0$, with $d(t)=\delta$ will tend to infinity. Hence $\Omega$ cannot be $\delta$-RUAS. 
\end{exmp}

Theorem \ref{thm:main} can be obtained as a corollary of Theorem \ref{thm:stable} and Theorem \ref{thm:lyap}. 

\begin{IEEEproof}[Proof of Theorem \ref{thm:main}]
By Theorem \ref{thm:stable}, $\Omega$ is $\delta'$-RUAS for any $\delta'\in[0,\delta)$. By Theorem \ref{thm:lyap}, there exists a neighborhood $D$ of $\Omega$ and a smooth $V:\,D\ra\Real$ such that $V$ satisfies conditions (1) and (2) in Definition \ref{def:lyap}.

Let
\begin{equation}\label{eq:bc}
B_c(x) = c- V(x),
\end{equation}
where $c>0$ is a scalar chosen sufficiently small such that $B_c(x)\ge 0$ implies $x\not\in U$. For instance, one can take $c$ be the maximum value of $V(x)$ in a compact neighborhood of $\Omega$ that does not overlap with $\overline{U}$. Then $B_c(x)$ also verifies all the conditions of a $\delta'$-robust barrier function. In particular, we have
\begin{align}
    \nabla B_c(x)\cdot(f(x)+d) & = -\nabla V(x)\cdot(f(x)+d) \notag\\
    & \ge \alpha_3(\alpha_2^{-1}(V(x))) \notag\\
    & = \alpha_3(\alpha_2^{-1}(c-B_c(x))) \notag\\
    & = \alpha_3(\alpha_2^{-1}(c)) > 0 \label{eq:pb}
\end{align}
for all $x$ such that $B_c(x)=0$ and all $d\in\delta\B$. 
\end{IEEEproof}

\begin{rem}
The construction of a Barrier function via a Lyapunov function is inspired by the work \cite{xu2015robustness} (see also \cite{ames2016control}) and \cite{ratschan2018converse}. In \cite{xu2015robustness}, the authors showed that if there exists of a barrier function $B:\,\Real^n\ra\Real$ satisfying the condition 
\begin{equation}\label{eq:b0}
\nabla B(x)\cdot(f(x))\ge -\alpha(B(x)),\quad \forall x\in D,
\end{equation}
for some open set $D$ containing $C=\set{x\in\Real^n:\,B(x)\ge 0}$ and extended class $\K$ function\footnote{A function $\alpha:(-b,a)\ra\Real$, $a,b>0$, is said to belong to extended class $\K$ if $\alpha$ is strictly increasing and $\alpha(0)=0$.} $\alpha$, then $C$ is asymptotically stable. This is straightforward to see because one can construct a Lyapunov function based on $B$ by $V(x)=0$ if $x\in C$ and $V(x)=-B(x)$ if $x\in D\backslash C$. The authors of \cite{xu2015robustness} also discussed robustness implied by condition (\ref{eq:b0}). The results of this section can be seen as a converse fact. We start with the assumption that a set $W$ is robustly safe and show that the closure of the robustly invariant reachable set $\Omega=\overline{R_{\delta}(W)}$ is   robustly asymptotically stable. Our proof of the latter fact is inspired by the work in \cite{ratschan2018converse}. A converse Lyapunov function is then used to construct a robust barrier function. 
\end{rem}

\begin{rem}\label{rem:barrier}
Condition (3) in Definition \ref{def:barrier} for a barrier function has different variants. The original work \cite{prajna2004safety} had a condition like (3) and the following variant \begin{equation}
    \label{eq:b1}
    \nabla B(x)\cdot f(x,d) \ge 0,\quad \forall (x,d)\in \mathcal{X}\times \mathcal{W},
\end{equation}
where $\mathcal{X}\times \mathcal{W}$ is the set on which $f$ is defined and $\mathcal{W}$ is an arbitrary disturbance set. According to \cite{prajna2004safety}, this variant makes the set of functions satisfying the barrier function conditions convex and amenable to computation by convex optimization. Condition (\ref{eq:b1}) appears to be restrictive (from a computational perspective) because it needs to be satisfied for all $(x,d)\in \mathcal{X}\times \mathcal{W}$. The authors of \cite{ames2016control} proposed (\ref{eq:b0}) as a variant. Following the construction $B(x)=-V(x)$ in the proof of Theorem \ref{thm:main}, we have
\begin{align}
    \nabla B(x)\cdot(f(x)+d) &= -\nabla V(x)\cdot(f(x)+d) \notag\\
    &\ge \alpha_3(\norm{x}_A) \notag\\
    &\ge \alpha_3(\alpha_2^{-1}(V(x))) \notag\\
    &= \alpha_3(\alpha_2^{-1}(-B(x))), \label{eq:nablaB}
\end{align}
for all $x\in D$ and $d\in \delta\mathbb{B}$. Defining $\alpha_0(s)=-\alpha_3(\alpha_2^{-1}(-s))$, we obtain
\begin{equation}
    \label{eq:b2}
    \nabla B(x)\cdot (f(x)+d) \ge -\alpha_0(B(x)),\quad \forall (x,d)\in D\times \delta'\mathbb{B}.
\end{equation}
While in the absence of disturbance (\ref{eq:b2}) appears in the same form as (\ref{eq:b0}), it has a subtle difference because $\alpha_0(s)$ in (\ref{eq:b2}) is not defined for $s>0$. Note that, since $B(x)=-V(x)$, $B(x)$ is never positive by this construction. Nonetheless, (\ref{eq:b2}) does match (\ref{eq:b0}) when $B(x)\le 0$ in the absence of disturbance. With the barrier function $B_c(x)$ defined in (\ref{eq:pb}) in the proof of Theorem \ref{thm:main}, we have
\begin{equation}
\label{eq:bc1}
\nabla B_c(x)\cdot(f(x)+d) \ge -\alpha_0(B_c(x))    
\end{equation}
with $\alpha_0(s)=-\alpha_3(\alpha_2^{-1}(c-s))$. Note that, compared with (\ref{eq:b2}),   $B_c(x)$ now can take positive value and $\alpha_0(s)$ is defined for $s\in(0,c]$ as well. Nonetheless, while  (\ref{eq:bc1}) agrees with (\ref{eq:b0}) for $B_c(x)\le 0$ in the absence of disturbance, it is in fact stronger than (\ref{eq:b0}) when $B_c(x)>0$ because $\alpha_0(s)<0$ for $s\in(0,c)$. This is not surprising because $B_c(x)$ is constructed using a Lyapunov function. Furthermore, in the absence of disturbance $d$, the strictly inequality (\ref{eq:pb}) established in the proof of Theorem \ref{thm:main} recovers condition (3) for the barrier function in Theorem \ref{thm:rats}. The author of \cite{ratschan2018converse} seems to be using this strict positiveness, as well as strict positiveness of $B$ on $W$, to indicate a robust barrier certificate. Here we formally define a robust barrier function by requiring condition (3) in Definition \ref{def:barrier} to hold with under disturbance. %
We also remark that, when the set $B(x)=0$ is compact, condition (3) in Theorem \ref{thm:rats} also holds under sufficiently small disturbance. The construction given by Theorem \ref{thm:main}, however, allows any disturbance of size $\delta'\in[0,\delta)$. 
\end{rem}

\begin{rem}\label{rem:reci}
Another commonly used class of barrier functions is called reciprocal barrier functions \cite{ames2016control}, inspired by barrier methods from optimization \cite{forsgren2002interior}. 
Given a set $C$ defined by 
$$
C= \set{x\in\Real^n:\,h(x)\ge 0},
$$
where $h:\,\Real^n\ra\Real$ is a continuously differentiable function, a reciprocal barrier function $B:\,C^\circ\ra\Real$, where $C^\circ= \set{x\in\Real^n:\,h(x)> 0}$ is the interior of $C$, such that
\begin{equation}
    \label{eq:rb1}
    \frac{1}{\alpha_1(h(x))} \le B(x) \le \frac{1}{\alpha_2(h(x))},
\end{equation}
\begin{equation}
    \label{eq:rb2}
    \nabla B(x)\cdot f(x) \le \alpha_3(h(x)),
\end{equation}
for all $x\in C^\circ$, where $\alpha_i$ ($i=1,2,3$) are class $\K$ functions. The reciprocal of the construction of barrier functions based on Lyapunov function directly gives a reciprocal barrier function. Let $h(x)=c-V(x)$ as in (\ref{eq:bc}) and $B(x)=\frac{1}{h}$. Then it is straightforward to verify that (\ref{eq:rb1}) is satisfied and (\ref{eq:rb2}) is robustly satisfied. 
\end{rem}

\section{Discrete-time converse barrier functions}
\label{sec:discrete}

Having built a link between Lyapunov functions and barrier functions, we extend the results in the previous section to the discrete-time setting and provide a converse theorem for discrete-time barrier function. The presentation parallels that of Section \ref{sec:main}, but formulated for discrete-time systems. We first present the preliminaries for discrete-time systems. 

\subsection{Preliminaries on discrete-time systems}

Consider a discrete-time dynamical system 
\begin{equation}\label{eq:dts}
    x(t+1) = f(x(t)),
\end{equation}
where $x(t)\in\Real^n$ for $t\in\mathbb{Z}^+\set{0,1,2,\cdots}$ and $f:\,\Real^n\ra\Real^n$ is assumed to be locally Lipschitz.

Given a scalar $\delta\ge 0$, a \emph{$\delta$-perturbation} of the dynamical system (\ref{eq:dts}) is described by the difference inclusion
\begin{equation}\label{eq:dtsp1}
    x(t+1) \in F_{\delta}(x(t)),
\end{equation}
where $F_{\delta}(x)=B_r(f(x))$, or equivalently 
\begin{equation}\label{eq:dtsp2}
    x(t+1) = f(x(t)) + d(t),
\end{equation}
where $d(t)\in\delta\B$ for each $t$. We denote system (\ref{eq:dts}) by $\dtS$ and its $\delta$-perturbation by $\dtSd$. Note that $\dtSd$ reduces to $\dtS$ when $\delta=0$. A solution of $\dtSd$ is a sequence denoted by $x(t;x_0,d)$ or $x(t)$, where $t=0,1,2,\cdots$ and $d(t)$ is a disturbance sequence. %

Since robustly safe sets, robustly invariant sets, and robust stability w.r.t. a closed set for $\dtS$ can be defined almost verbatim as for continuous-time systems, by replacing solutions of $\Sd$ with that of $\dtSd$, they are omitted. We define discrete-time barrier and Lyapunov functions as follows. 

\begin{defn}[Discrete-time robust barrier function]\label{def:dtbarrier}
Given sets $W,U\subset\Real^n$, a continuously differentiable function  $B:\,\Real^n\ra\Real$ is said to be a \textit{$\delta$-robust barrier function} for $W$ and $U$ if the following conditions are satisfied: 
\begin{enumerate}[(1)]
    \item $B(x)\ge 0$ for all $x\in W$;
    \item $B(x)<0$ for all $U$; and 
    \item $B(f(x)+d)\ge 0$ for all $x$ such that $B(x)\ge 0$ and all $d\in\delta\B$.
\end{enumerate}
\end{defn}

\begin{prop}[Sufficiency of discrete-time barrier functions]
If there exists a $\delta$-robust barrier function for $(W,U)$, then $W$ is $\delta$-robustly safe w.r.t. $U$. 
\end{prop}

\begin{IEEEproof}
The conclusion follows from the fact that the set $C=\set{x\in\Real^n:\,B(x)\ge 0}$ is $\delta$-robustly invariant and $C\cap U\neq\emptyset$. 
\end{IEEEproof}

\begin{rem}
Condition (3) in Definition \ref{def:dtbarrier} for a discrete-time barrier function appears to be weaker than the ones used in practice. For instance, the following condition was proposed in \cite{ahmadi2019safe}:
\begin{equation}\label{eq:dtb0}
    B(f(x)) - B(x) \ge -\alpha(B(x)), \quad x\in D,
\end{equation}
where $D\supseteq C=\set{x\in\Real^n:\,B(x)\ge 0}$  and $\alpha$ is class $\K$ function satisfying $\alpha(r)<r$ when $r>0$. Note that one needs to extend the definition of $\alpha$ to $(-b,0)$ for some $b>0$ if the set $D$ contains $x$ such that $B(x)<0$. A special case of (\ref{eq:dtb0}) is given by $\alpha(r)=\eta r$ for $\eta\in(0,1)$. When $\eta=1$ and $D=C$, we obtain condition (3) of Definition \ref{def:dtbarrier}. When $\eta=0$, we obtain a condition that is stronger than (\ref{eq:dtb0}) on $C$:
\begin{equation}\label{eq:dtb1}
    B(f(x)) - B(x) \ge 0, \quad x\in C,
\end{equation}
which clearly implies (\ref{eq:dtb0}) for any $\alpha\in\K$ and $D=C$ because $\alpha(B(x))\ge 0$ for $B(x)\ge 0$. Similar to Remark \ref{rem:barrier} on continuous-time barrier functions, the construction of discrete-time converse barrier functions by Lyapunov functions below in fact satisfy an even stronger form  \begin{equation}\label{eq:dtb2}
    B(f(x)) - B(x) \ge -\alpha_0(B(x)), \quad x\in D,
\end{equation}
where $D$ is an open neighborhood of $C$ and $\alpha_0(s)\le 0$ for all $s$ in its domain. Clearly, (\ref{eq:dtb2}) implies both (\ref{eq:dtb1}) and (\ref{eq:dtb0}). 
\end{rem}

\begin{defn}\label{def:dtlyap}
Let $D\subset\Real^n$ be an open set containing a closed set $A\subset\Real^n$. A continuously differentiable function $V:\,D\ra\Real$ is said to be a \textit{$\delta$-robust Lyapunov function} for $\dtS$ w.r.t. $A$ if the following two conditions are satisfied:
\begin{enumerate}[(1)]
    \item there exist class $\K$ functions $\alpha_1$ and $\alpha_2$ such that 
    $$
    \alpha_1(\norm{x}_A) \le V(x) \le \alpha_2(\norm{x}_A)
    $$
    for all $x\in D$; and 
    \item there exists a class $\K$ functions $\alpha_3$ such that
    $$
    V(f(x)+d) - V(x)\le -\alpha_3(\norm{x}_A)
    $$
    for all $x\in D$ and $d\in \delta\B$. 
\end{enumerate}
\end{defn}

There are also Lyapunov characterizations of set stability for discrete-time systems. 

\begin{thm}[Lyapunov characterization of set\label{thm:lyapdt} stability for $\dtS$ \cite{jiang2002converse}]
A closed set $A\subset\Real^n$ is $\delta$-RUAS for $\dtS$ if and only if there exists a $\delta$-robust Lyapunov function for $\dtS$ w.r.t. $A$. 
\end{thm}

\subsection{Converse barrier functions via Lyapunov functions for discrete-time systems}

The notation and definitions for reachable sets remain the same, with continuous-time solutions replaced with discrete-time ones. We define $R_\delta^t(x_0)$, $R_\delta(x_0)$, $R_{\delta}(W)$ and $\Omega=\overline{R_\delta(W)}$ as in Section \ref{sec:main}, replacing continuous-time solutions with discrete-time ones. The following is a discrete-time version of Lemma \ref{lem:control}. %

\begin{lem}\label{lem:controldt}
Fix any $\delta'\in(0,\delta)$. Let $K\subset\Real^n$ be a compact set. Then there exists some $r=r(K,\delta',\delta)>0$ such that the following holds: if there is a solution $x$ of $\S_{\delta'}$ such that $x(s)\in K$ for all $s\in[0,T]$, where $T\ge 1$, then for any $y_0\in B_r(x(0))$ and any $y_1\in B_r(x(T))$, we have $y_1\in R_{\delta}^T(y_0)$, i.e., $y_1$ is reachable at $T$ from $y_0$ by a solution of $\Sd$. Furthermore, if $f$ is globally Lipschitz, $r$ can be chosen to be independent of $K$. 
\end{lem}

\begin{IEEEproof}
Let 
$$
y(s) = x(s) + \frac{s}{T}[y_1 - x_1 + (x_0-y_0)] + (y_0-x_0),
$$
for $s\in\set{0,1,\cdots,T}.$ Then $y(0)=y_0$ and $y(T)=y_1$. Furthermore, 
\begin{align*}
    \norm{y(s)-x(s)} &\le \norm{y_1-x_1}\frac{s}{T} + \norm{y_0-x_0}(1-\frac{s}{T}) \\ 
    &\le r (\frac{s}{T} + 1  - \frac{s}{T}) = r,
\end{align*}
for all $s\in\set{0,1,\cdots,T}.$ Hence
\begin{align*}
    &\norm{y(s+1)-f(y(s))} \\
    &= \norm{y(s+1)-x(s+1)} + \norm{x(s+1) - f(x(s))}\\
    &\qquad + \norm{f(x(s)) - f(y(s))} \\
    & \le \frac{r}{T} + \delta' + L r,\quad s\in \set{0,1,\cdots,T-1},
\end{align*}
where we used the triangle inequality, the fact that $x$ is a solution of $\S_{\delta'}$, and Lipschitz continuity of $f$ on the set $B_r(K)$. By picking $r$ sufficiently small such that $r + \delta' + L r<\delta$, then we have
$\norm{y(s+1)-f(y(s))}<\delta$ for all $s\in \set{0,1,\cdots,T-1}$. Thus $y$ is a solution of $\dtSd$ and the conclusion follows. Note that the choice of $r$ only depends on $K$, $\delta'$, and $\delta$. If $f$ is globally Lipschitz, the dependence on $K$ can be removed. 
\end{IEEEproof}

The following is a discrete-time version of Theorem \ref{thm:stable}. 

\begin{thm}[Robustly invariant sets are robustly asymptotically stable]\label{thm:stabledt}
If Assumption \ref{as:main} holds, then $\Omega$ is $\delta'$-RUAS for $\dtS$ for any $\delta'\in[0,\delta)$.
\end{thm}

\begin{IEEEproof}
(1) For any $\eps>0$, let $r$ be from Lemma \ref{lem:controldt}, applied to the set $B_{\eps}(\Omega)$ and scalars $\delta'$ and $\delta$. Pick $\delta_\eps=r$. Let $x$ be any solution of $\dtS_{\delta'}$ such that $\norm{x(0)}_\Omega<\delta_\eps$. We show that $\norm{x(t)}_\Omega<\eps$ for all $t\ge 0$. Suppose that this is not the case. Then $\norm{x(k)}_\Omega\ge\eps$ for some $k\ge 1$.  
Since $\norm{x(0)}_\Omega<r$, we can always pick $y_0\in R_\delta(W)$ such that $y_0\in B_{r}(x(0))$ by the triangle inequality. By Lemma \ref{lem:controldt}, there exists a solution of $\dtSd$ from $y_0\in R_\delta(W)$ to $x(k)\not\in\Omega$. This contradicts that $R_\delta(W)$ is $\delta$-robustly invariant. 

(2) Fix any $\eps_0>0$. Following part (1), choose $\delta_{\eps_0}$ such that $\norm{x(0)}_\Omega<\delta_{\eps_0}$ implies $\norm{x(t)}_\Omega<\eps_0$ for any solution $x(t)$ of $\dtS_{\delta'}$. Let $r$ be chosen according to Lemma \ref{lem:controldt} with the set $B_{\eps_0}(\Omega)$ and scalars $\delta'$ and $\delta$. Choose $\rho\in(0,r)$. Let $x$ be any solution of $\dtS_{\delta'}$. We show that $\norm{x(0)}_A<\rho$ implies that $x(t)\in R_\delta(W)$ for all $t\ge 1$. Suppose that this is not the case. Then there exists some $k\ge 1$ such that $x(k)\in\partial \Omega$ or $x(k)\not\in\Omega$. In either case, we can pick $y_1\in B_r(x(k))$ such that $y_1\not\in\Omega$ and $y_0\in B_r(x(0))$ such that $y_0\in  R_\delta(W)$. By Lemma \ref{lem:controldt}, there exists a solution of $\dtSd$ from $y_0\in R_\delta(W)$ to $y_1\not\in\Omega$. This contradicts that $R_\delta(W)$ is $\delta$-robustly invariant. Hence $x(t)\in R_\delta(W)$ for all $t\ge 1$. This completes part (2) of the definition of $\delta'$-RUAS.  %
\end{IEEEproof}

\begin{thm}[Robustly safe sets admit robust discrete-time barrier functions]\label{thm:maindt}
Suppose that Assumption \ref{as:main} holds. If either $\Omega$ is compact or $f$ is globally Lipschitz, then for any $\delta'\in(0,\delta)$, there exists a $\delta'$-robust barrier function for $(W,U)$.
\end{thm}

\begin{IEEEproof}
By Theorem \ref{thm:stabledt}, $\Omega$ is $\delta'$-RUAS for $\dtS$ with any $\delta'\in[0,\delta)$. By Theorem \ref{thm:lyapdt}, there exists a neighborhood $D$ of $\Omega$ and a smooth $V:\,D\ra\Real$ such that 
$$
\alpha_1(\norm{x}_\Omega) \le V(x) \le \alpha_2(\norm{x}_\Omega),
$$
and
$$
V(f(x)+d) - V(x) \le -\alpha_3(\norm{x}_\Omega),
$$
for all $x\in D$ and $d\in\delta' \mathbb{B}$, where $\alpha_i$ ($i=1,2,3$) are class $\K$ functions. 
Define 
$$
B(x) = -V(x),\quad x\in D. 
$$
It is straightforward to verify that $B$ satisfies conditions (1)--(3) of Definition \ref{def:dtbarrier} for a $\delta'$-robust discrete-time barrier function. 
\end{IEEEproof}

\begin{rem}
By the construction of the barrier function $B(x)$ in the proof of Theorem \ref{thm:maindt}, we in fact have a stronger condition than condition (3) in Definition \ref{def:dtbarrier}: 
\begin{equation}
    \label{eq:barrierdt}
    B(f(x)+d) - B(x) \ge -\alpha(B(x)),
\end{equation}
for all $x\in D$ and $d\in\delta'\mathbb{B}$, where $\alpha$ is defined and increasing on $(-a,0]$ for some $a>0$ and $\alpha(0)=0$. 
\end{rem}

\begin{rem}
Similar to Remark \ref{rem:reci}, we can construct discrete-time reciprocal barrier functions via Lyapunov functions. A discrete-time reciprocal barrier function \cite{agrawal2017discrete} $B:\,C^\circ\ra\Real$ satisfies 
\begin{equation}
    \label{eq:dtrb1}
    \frac{1}{\alpha_1(h(x))} \le B(x) \le \frac{1}{\alpha_2(h(x))},
\end{equation}
\begin{equation}
    \label{eq:dtrb2}
     B(f(x)) - B(x)\le \alpha_3(h(x)),
\end{equation}
for all $x\in C^\circ$, where $C^\circ$ is the interior of the set 
$
C= \set{x\in\Real^n:\,h(x)\ge 0}
$
for some continuously differentiable function $h:\,\Real^n\ra\Real$. Clearly, $h(x)=c-V(x)$ for some sufficiently small $c>0$ and $B(x)=\frac{1}{h}$ satisfy the above conditions robustly.   
\end{rem}

\section{Conclusions}
\label{sec:conclusion}

In this paper, we established a connection between Lyapunov functions and barrier functions. We proved that for all robustly safe dynamical systems, the closure of the robust reachable set of the robustly safe set must be robustly asymptotically stable. The converse Lyapunov function theory can then be brought to bear to yield a robust barrier function. We made remarks on several variants of the barrier function conditions and showed that they can all be satisfied by the construction of barrier functions using Lyapunov functions. We also formulated the results for discrete-time in a similar fashion. 

For future work, it would be interesting to investigate how the viewpoint of robust barrier functions via Lyapunov functions can be utilized in practice. Potentially all the computational techniques for searching Lyapunov functions can be used for searching barrier functions. The key technical challenge, however, seems to be that, while safety requirements can be specified rather arbitrarily by a designer (e.g., by defining the unsafe region $U$ and safe initial region $W$ in this paper), the barrier function conditions are only met at the boundary of reachable set from the safe initial region $W$, if this set $W$ can indeed be certified to be safe. While the computing of reachable sets can be highly nontrivial, it would be interesting to investigate whether the adaptive refinement techniques for computing maximal controlled invariant sets (see, e.g., \cite{li2017invariance}), combined with computational techniques for constructing barrier functions (see, e.g., \cite{djaballah2017construction}), can be used to determine a smaller set on which (control) barrier functions can be algorithmically constructed. A related theoretical question is that whether such procedures can be approximately complete in the sense that any $\delta$-robustly safe sets admit a computable $\delta'$-robust barrier certifications for any $\delta'\in[0,\delta)$. In view of the results of this paper, such questions can hopefully be answered in a unified fashion in regard to Lyapunov functions for set stability and barrier functions for safety. 

\section{Acknowledgments}

The author would like to thank Stefan Ratschan for very helpful comments and, in particular, for pointing out the need for the strict inequality condition (3) in Definition \ref{def:barrier}.  
 
\bibliographystyle{IEEEtran}  
\bibliography{references}  %

\begin{thebibliography}{10}
\providecommand{\url}[1]{#1}
\csname url@samestyle\endcsname
\providecommand{\newblock}{\relax}
\providecommand{\bibinfo}[2]{#2}
\providecommand{\BIBentrySTDinterwordspacing}{\spaceskip=0pt\relax}
\providecommand{\BIBentryALTinterwordstretchfactor}{4}
\providecommand{\BIBentryALTinterwordspacing}{\spaceskip=\fontdimen2\font plus
\BIBentryALTinterwordstretchfactor\fontdimen3\font minus
  \fontdimen4\font\relax}
\providecommand{\BIBforeignlanguage}[2]{{%
\expandafter\ifx\csname l@#1\endcsname\relax
\typeout{** WARNING: IEEEtran.bst: No hyphenation pattern has been}%
\typeout{** loaded for the language `#1'. Using the pattern for}%
\typeout{** the default language instead.}%
\else
\language=\csname l@#1\endcsname
\fi
#2}}
\providecommand{\BIBdecl}{\relax}
\BIBdecl

\bibitem{prajna2007convex}
S.~Prajna and A.~Rantzer, ``Convex programs for temporal verification of
  nonlinear dynamical systems,'' \emph{SIAM Journal on Control and
  Optimization}, vol.~46, no.~3, pp. 999--1021, 2007.

\bibitem{prajna2005necessity}
------, ``On the necessity of barrier certificates,'' \emph{IFAC Proceedings
  Volumes}, vol.~38, no.~1, pp. 526--531, 2005.

\bibitem{prajna2004safety}
S.~Prajna and A.~Jadbabaie, ``Safety verification of hybrid systems using
  barrier certificates,'' in \emph{International Workshop on Hybrid Systems:
  Computation and Control}.\hskip 1em plus 0.5em minus 0.4em\relax Springer,
  2004, pp. 477--492.

\bibitem{prajna2007framework}
S.~Prajna, A.~Jadbabaie, and G.~J. Pappas, ``A framework for worst-case and
  stochastic safety verification using barrier certificates,'' \emph{IEEE
  Transactions on Automatic Control}, vol.~52, no.~8, pp. 1415--1428, 2007.

\bibitem{wieland2007constructive}
P.~Wieland and F.~Allg{\"o}wer, ``Constructive safety using control barrier
  functions,'' \emph{IFAC Proceedings Volumes}, vol.~40, no.~12, pp. 462--467,
  2007.

\bibitem{agrawal2017discrete}
A.~Agrawal and K.~Sreenath, ``Discrete control barrier functions for
  safety-critical control of discrete systems with application to bipedal robot
  navigation,'' in \emph{Robotics: Science and Systems}, 2017.

\bibitem{cheng2019end}
R.~Cheng, G.~Orosz, R.~M. Murray, and J.~W. Burdick, ``End-to-end safe
  reinforcement learning through barrier functions for safety-critical
  continuous control tasks,'' in \emph{Proceedings of the AAAI Conference on
  Artificial Intelligence}, vol.~33, 2019, pp. 3387--3395.

\bibitem{ahmadi2019safe}
M.~Ahmadi, A.~Singletary, J.~W. Burdick, and A.~D. Ames, ``Safe policy
  synthesis in multi-agent pomdps via discrete-time barrier functions,'' in
  \emph{2019 IEEE 58th Conference on Decision and Control (CDC)}.\hskip 1em
  plus 0.5em minus 0.4em\relax IEEE, 2019, pp. 4797--4803.

\bibitem{ames2016control}
A.~D. Ames, X.~Xu, J.~W. Grizzle, and P.~Tabuada, ``Control barrier function
  based quadratic programs for safety critical systems,'' \emph{IEEE
  Transactions on Automatic Control}, vol.~62, no.~8, pp. 3861--3876, 2017.

\bibitem{xu2015robustness}
X.~Xu, P.~Tabuada, J.~W. Grizzle, and A.~D. Ames, ``Robustness of control
  barrier functions for safety critical control,'' \emph{IFAC-PapersOnLine},
  vol.~48, no.~27, pp. 54--61, 2015.

\bibitem{tee2009barrier}
K.~P. Tee, S.~S. Ge, and E.~H. Tay, ``Barrier lyapunov functions for the
  control of output-constrained nonlinear systems,'' \emph{Automatica},
  vol.~45, no.~4, pp. 918--927, 2009.

\bibitem{romdlony2016stabilization}
M.~Z. Romdlony and B.~Jayawardhana, ``Stabilization with guaranteed safety
  using control {Lyapunov}--barrier function,'' \emph{Automatica}, vol.~66, pp.
  39--47, 2016.

\bibitem{wisniewski2016converse}
R.~Wisniewski and C.~Sloth, ``Converse barrier certificate theorems,''
  \emph{IEEE Transactions on Automatic Control}, vol.~61, no.~5, pp.
  1356--1361, 2016.

\bibitem{ratschan2018converse}
S.~Ratschan, ``Converse theorems for safety and barrier certificates,''
  \emph{IEEE Transactions on Automatic Control}, vol.~63, no.~8, pp.
  2628--2632, 2018.

\bibitem{khalil2002nonlinear}
H.~K. Khalil, \emph{Nonlinear Systems}.\hskip 1em plus 0.5em minus 0.4em\relax
  Prentice-Hall, 2002.

\bibitem{taly2009deductive}
A.~Taly and A.~Tiwari, ``Deductive verification of continuous dynamical
  systems,'' in \emph{IARCS Annual Conference on Foundations of Software
  Technology and Theoretical Computer Science}.\hskip 1em plus 0.5em minus
  0.4em\relax Schloss Dagstuhl-Leibniz-Zentrum f{\"u}r Informatik, 2009.

\bibitem{gard1980strongly}
T.~C. Gard, ``Strongly flow-invariant sets,'' \emph{Applicable Analysis},
  vol.~10, no.~4, pp. 285--293, 1980.

\bibitem{blanchini1999set}
F.~Blanchini, ``Set invariance in control,'' \emph{Automatica}, vol.~35,
  no.~11, pp. 1747--1767, 1999.

\bibitem{wilson1969smoothing}
F.~W. Wilson, ``Smoothing derivatives of functions and applications,''
  \emph{Transactions of the American Mathematical Society}, vol. 139, pp.
  413--428, 1969.

\bibitem{lin1996smooth}
Y.~Lin, E.~D. Sontag, and Y.~Wang, ``A smooth converse {Lyapunov} theorem for
  robust stability,'' \emph{SIAM Journal on Control and Optimization}, vol.~34,
  no.~1, pp. 124--160, 1996.

\bibitem{liu2019robust}
J.~Liu, ``Robust decidability of sampled-data control of nonlinear systems with
  temporal logic specifications,'' \emph{arXiv preprint arXiv:1903.06368},
  2019.

\bibitem{forsgren2002interior}
A.~Forsgren, P.~E. Gill, and M.~H. Wright, ``Interior methods for nonlinear
  optimization,'' \emph{SIAM review}, vol.~44, no.~4, pp. 525--597, 2002.

\bibitem{jiang2002converse}
Z.-P. Jiang and Y.~Wang, ``A converse lyapunov theorem for discrete-time
  systems with disturbances,'' \emph{Systems \& Control Letters}, vol.~45,
  no.~1, pp. 49--58, 2002.

\bibitem{li2017invariance}
Y.~Li and J.~Liu, ``Invariance control synthesis for switched nonlinear
  systems: An interval analysis approach,'' \emph{IEEE Transactions on
  Automatic Control}, vol.~63, no.~7, pp. 2206--2211, 2018.

\bibitem{djaballah2017construction}
A.~Djaballah, A.~Chapoutot, M.~Kieffer, and O.~Bouissou, ``Construction of
  parametric barrier functions for dynamical systems using interval analysis,''
  \emph{Automatica}, vol.~78, pp. 287--296, 2017.

\end{thebibliography}

\end{document}